\documentclass{article}
\usepackage[utf8]{inputenc}
\usepackage[centertags]{amsmath}
\usepackage{amsfonts,amssymb,amsthm}

\title{Kurt G\"{o}del and the Logic of Concepts}

\author{
Jovana Kosti\'{c}\footnote{Supported by the Ministry of Education, Science and Technological Development of the Republic of Serbia through the financing of the scientific research on the University of Belgrade - Faculty of Philosophy (contract number 451-03-47/2023-01/ 200163).} \\ {\small Faculty of Philosophy, University of Belgrade} \\[-2mm] {\small email: kostic.jovana345@gmail.com}
\and
Slobodan Vujo\v{s}evi\'{c} \\{\small Mathematical Institute, SASA,  Belgrade} \\ [-2mm] {\small email: vujosevic21@gmail.com}
}
\date{}

\begin{document}

\maketitle

\begin{abstract}
The literature dealing with G\"{o}del's legacy is largely preoccupied with challenging his philosophical views, regarding them as outdated. We believe that such an approach prevents us from seeing G\"{o}del's views in the right light and understanding their rationale. In this article, his views are discussed in the philosophical realm in which he himself understood them. We explore the consequences of G\"{o}del's incompleteness theorems for the question of the objectivity of mathematics and its epistemology. Taking set theory as the paradigm of formal mathematical theories, we examine the relationship between its incompleteness and extensionality. We argue, based on his philosophical views, that G\"{o}del believed incompleteness can be overcome only by some intensional considerations about concepts from the basis of mathematical theories. These considerations should eventually lead to founding the so-called logic of concepts.
\end{abstract}

\section{Introduction}  

By the end of the XIX century, mathematics was fully extensionalized. All mathematical concepts were taken to be reducible to set-theoretic ones. Logic adjusted to the needs of founding so-conceived mathematics. The properties of concepts and their logical structure were practically excluded from the study, contrary to G\"{o}del's view that ``logic is a theory of concepts'' (Kurt G\"{o}del, 1906 - 1978). The exceptions were made by sporadic attempts at founding an {\it intensional logic}, a philosophical discipline that examines the relationship between a name and its denotation, and deals with related problems, such as the epistemic value of the statement ``Morning star is Evening star'', or the truth value of the statement ``The present king of France is bald''. G\"{o}del had a very different idea of how intensional logic should look like. He imagined it as a theory that deals with the formal properties of concepts and clarifies their relation to sets. For this reason, we call the theory envisioned by G\"{o}del {\it the logic of concepts}, to distinguish it from the theories that nowadays fall under the name of intensional logic.
    
As far as G\"{o}del's legacy allows, we will here analyze the reasons that led him to the idea of founding the logic of concepts. We believe the decisive reason was G\"{o}del's dissatisfaction with the current state of the fundamental mathematical theory --  set theory. As an incomplete theory, it eventually assumed the form of an ``if-then science''. It has ramified into multiple theories exploring different ways of redeeming its incompleteness. This is at odds with G\"{o}del's belief that the mathematical world is objective and unique and that knowledge of it can only be extended from verified truths rather than unverified hypotheses. In contrast to empirical knowledge, which is based on hypotheses that accord with observations, the starting point in mathematics is indubitable truth. 

G\"{o}del might have thought that the logic of concepts, as an intensional theory, will not be limited by his incompleteness theorems and forced to rely on unverified hypotheses. It might be even possible for it to prove its consistency. This would be a ``mathematically rigorous'' and invincible argument in favor of G\"{o}del's principal belief in the objective existence of the world of concepts.

\section{Incompleteness theorems and Platonism}

Among the most intriguing philosophical enigmas of our time, whose consequences are still not thoroughly investigated after almost a century, are the theorems stating the incompleteness of the most important mathematical formal systems. Many philosophers and logicians have explored the significance of these theorems. However, they would often find a reason to distance themselves from the interpretation their author, Kurt G\"{o}del, ascribed them. G\"{o}del believed that his results provide a strong argument for {\it the objective existence of a rationally organized world of concepts, which can be to some extent described by a deductive system but cannot be changed or manipulated} (cf. \cite{godel3}, p. 320). Despite its authenticity, this interpretation of G\"{o}del's results was not acceptable to {\it the spirit of his time} and is hardly acceptable even today. It is usually looked upon ``with materialistic and positivistic prejudice'' as old-fashioned and obsolete.  

G\"{o}del's results have attracted huge attention, they lifted logic to a higher level and led to the establishment of new areas of study, which created possibilities for thousands of mathematicians, logicians, and philosophers. But, with regard to their interpretation, G\"{o}del was alone. He was not a philosopher, but a mathematician and logician who faced important philosophical questions working in his fields. He needed help from philosophers to make sense of them, and as he was not met with understanding, he turned to philosophy himself. Due to the occasion of the Second World War and possibly his introvert personality, G\"{o}del hasn't built a standard university career. He hasn't gathered a group of students that would be able to continue his work and develop it consistently and systematically. The metaphysical essence of G\"{o}del's views was foreign to the naturalistic aspirations of the science and philosophy of the time. This might be why G\"{o}del withdrew from academia very early. He gave his last lecture already in 1951. His subsequent notes on philosophical questions were intended to be kept private. They do not represent discussions with philosophers, but the development of G\"{o}del's philosophical standpoint built on his mathematical results.

G\"{o}del explicitly says that he formed an opinion on the existence of mathematical objects and the objectivity of mathematical propositions already during his studies (see \cite{godel4}, p. 444). At that time, he became interested in philosophy, especially in Kant with whom he shared epistemic optimism and the belief that mathematical knowledge can be indefinitely extended. The philosopher that made the biggest impact on G\"{o}del was Leibniz. Under his influence, G\"{o}del adopted a strong version of rationalism that includes the belief that philosophy can be established as a deductive system, following the example of mathematics and logic. Some authors argue that G\"{o}del's philosophical views were changeable and that at first, they manifested a moderate form of Platonism (see \cite{goldfarb}). This refers, in the first place, to the authors that prepared G\"{o}del's {\it Nachlass} for publication. They pushed G\"{o}del's philosophical ideas in the background, because of their metaphysical nature, and made it appear as if they did not have any significant influence on his logical and mathematical results. 

However, G\"{o}del was not a logician and mathematician in the ordinary sense; his work in these sciences was driven by a clear and well-worked-out vision grounded in the best philosophical tradition. In the preparation of his collected works, this vision was systematically neglected. It is presented as a side issue with no important role in founding and motivating his results. On the other hand, G\"{o}del's mathematical opus suggests that his program of mathematical investigations is made to implement and confirm his philosophical views. G\"{o}del's mathematics is a result of his philosophical assumptions, it follows from them and, hence, agrees with them. After Leibniz, there is no other mathematician or philosopher whose work has this quality. This applies to his incompleteness theorems as well. It was not by chance that G\"{o}del chose arithmetic, as the simplest mathematical theory, to examine the relationship between the concepts of {\it truth} and {\it provability}. This gives us strong reasons to think that G\"{o}del's Platonistic vision of mathematics was not changing in essence. What was changing were the circumstances in which he presented and defended his views, in which Platonism always had a central place -- till the 40s as {\it mathematical realism}, and later as {\it conceptual realism}. 

As will be explained later on, incompleteness theorems make for a strong argument for the objectivity of mathematical truth. As G\"{o}del remarked, it is the belief in this objectivity that motivated and even made possible his most important logical results -- the completeness theorem for predicate logic, as well as the incompleteness theorems for arithmetic and stronger theories. These results show that G\"{o}del deemed the problem of completeness of logical and mathematical formal systems very important. This is stressed in the lecture notes for the elementary logic course he gave in 1939 (see \cite{miloskosta}). The fact that sufficiently strong mathematical formal systems are bound to be incomplete, must have been regarded by G\"{o}del as revealing about the nature of mathematical knowledge. It is true that G\"{o}del at first expresses his views only by insisting on the difference between the concepts of {\em truth} and {\em proof} as the consequence of his theorems. But it should be remembered that this happens in the presence of philosophers of the Vienna Circle, of whose positivistic skepticism towards concepts with metaphysical content G\"{o}del was well aware. For them, no concept of truth transcending that of provability would be allowed. Otherwise, ``metaphysics might enter through the backdoor the edifice of their philosophy from which it was expelled through the main door'' (see \cite{karnap2}). 

Being aware of the outright hostility towards his metaphysical conception, relying solely on the nature of mathematics akin to it, G\"{o}del expresses his views modestly and carefully. At the same time, he seeks to establish some kind of {\it conclusive argument for the objectivity of mathematics}. In what follows, we present an overview of these G\"{o}del's attempts that include investigations of the philosophical consequences of his theorems for set theory, philosophy of mathematics, epistemology, and the logic of concepts.

\section{Objectivity of set theory}

According to G\"{o}del, the objectivity of mathematical propositions has to do with their analytic nature and close connection to physical reality. ``The {\it number} of petals is just as objective as the {\it color} of a flower'' (\cite{wang1}, 9.2.41). Similarly, the analytic form of law according to which objects move is equally objective as the material of which they are made. We arrive at it by abstracting from all the physical parameters, including those of space and time. By basing the objectivity of mathematical propositions on the intuition of sets (and, later, objectively existing concepts), which is free from any intuition of space or time, G\"{o}del departs from Kant's position in the philosophy of mathematics. In Kant's view, mathematics has to do with the way we build knowledge of the physical world, and is thus closely related to epistemology. G\"{o}del's departure from this standpoint makes mathematics closer to metaphysics and represents the revival of some prior philosophers, first of all, Leibniz.

Till the end of the XIX century, mathematical propositions were taken to have arithmetic or geometric content. But almost suddenly this has changed. Mathematical objects are reduced to sets, and propositions about them are given set-theoretic interpretation. Geometrical or arithmetical intuition of concepts is replaced by the intuition of sets. Mathematics has become an {\it extensional} science. Concepts remained in the background of so-conceived mathematics forming its metatheory, but the intention of ``objectifying'' them in sets remained present. This is not to say that till this moment mathematics was {\it intensional} in the sense of dealing with concepts and their properties, and it can hardly be said that it was transformed from conceptual into set-theoretic science. But some intensional aspects of the basic mathematical concepts that used to be taken into account, started to be neglected or denied. 

The tendency of understanding mathematical concepts as sets existed long before Cantor, but his precise treatment of multiplicities as unities pushed forward this idea. Sets are, undoubtedly, the most rewarding and illuminating way of understanding concepts. They provide perhaps the only link between the world of concepts and the physical world (cf. \cite{wang1}, 8.2.4). This change in understanding of mathematical content led to the change in {\it mathematical epistemology}. The main cognitive tools in mathematics became {\it definability}, {\it provability}, and {\it computability}. Thanks to that, mathematics achieved the precision and exactness we are now so used to. The concept of set was soon submitted to the same treatment, which resulted in establishing the formal set theory, better known as the Zermelo-Fraenkel set theory ($ZF$ or $ZFC$ if it includes the Axioms of Choice).

\subsection{Iterative conception of sets}

G\"{o}del's understanding of the set universe is tied to the {\it iterative conception of sets} (see \cite{godel3}, pp. 45-54). This conception is based on Cantor's understanding of a set as ``any collection into a whole of definite, well-distinguished objects of our intuition or thought'' (translation taken from \cite{bulos}, p. 215). This idea was employed multiple times (by Zermelo and von Neumann), but with G\"{o}del, it has become an indispensable part of the work in set theory. His first formulation of the iterative conception G\"{o}del named {\it the simple theory of types}. In contrast to Russell's theory, it allows the ``mixing'' of types, i.e., the existence of sets containing sets from different lower types. G\"{o}del presented the construction of the set-theoretic universe in a lecture delivered in 1933, during his first visit to America (see the manuscript of the lecture: \cite{godel3}, pp. 45-54).

The iterative conception of sets is of great philosophical interest. It offers a description of sets that can be understood objectively, as G\"{o}del does, but can also be accommodated to other philosophies of mathematics. Starting from the empty set, a transfinite iteration of the power set operation produces the family of sets $(V_{\alpha}, \alpha \in Ord)$ whose members are: $V_0 = \emptyset; V_{\alpha + 1} = P(V_\alpha)$, and $V_\alpha = \bigcup_{\beta < \alpha} V_\beta$, for a limit ordinal $\alpha$. The union of this family, denoted by $V$, represents {\it the universe of sets}. This universe $V$ is not a set, but a class. Still, it makes an appropriate model of the set-theoretic axioms because the meaning of the formula $x \in V$, which is that {\it $x$ is a set}, can also be expressed by the formula $\exists \alpha ( x \in V_\alpha)$, which is a formula of $ZF$ theory. Thus the statements about the properties of sets in the set universe need not be taken to describe the objective world, but can also be formally understood. 

The infinite ordinals used in the description of the levels of $V$ are easily definable in its formal construction provided by $ZF$ theory. However, if we aim to describe the objectively existing universe that justifies the axioms of the theory without being previously given the class of ordinals, then its description may appear circular. To be able to formulate the axioms for $V_\alpha$, we use the notion of an ordinal $\alpha$. G\"{o}del solves this problem in a typical Platonist vein, by freely using very strong mathematical arguments including the second-order Axiom of Choice and a variant of {\it Reflection Principle} saying that for every well-defined property of the set universe, there is a set with this property (see \cite{milos}, p. 203-240). 

The soundness of G\"{o}del's informal construction of the set universe is corroborated by an iterative construction published by George Boolos in 1971 (see \cite{bulos}). He offers a formal description of the stages of the iterative formation of sets - {\it the stage theory}. The language of this theory is strictly separated from the language of set theory. Besides the axioms describing the properties of the stages in iterative construction, the theory contains the second-order axiom of complete induction. It validates all $ZF$ axioms and shows that this cannot be done for the Axiom of Choice. The manuscript of the lecture in which G\"{o}del presents his simple type theory reveals his concern for the justification of the Axiom of Choice. It is not clear that Boolos' construction resolves it (see \cite{milos}). 

{\it Contrary to formalists that use the formulas of $ZF$ theory to define the set universe, G\"{o}del starts with the set universe built according to the iterative conception of sets and uses it to justify the $ZF$ axioms}. This difference is of crucial importance for understanding his objectivist conception of mathematics. According to it, the axioms of set theory are supposed to express the truth about the existing mathematical reality, i.e. they need to be confirmed, and not only assumed, as in the formalists' view. On the other hand, a large majority of mathematicians deny that mathematics is a system of truths about the world of mathematical objects. They accept that this holds only for its smaller part (according to their temperament, as G\"{o}del saw it), and understand the rest of mathematics in a hypothetical sense: some conclusions can be derived from the given assumptions whose justification ultimately depends on their applicability. In this way, mathematics is reduced to an empirical science.

\subsection{Justification of set axioms}

According to G\"{o}del, every mathematical theory needs to satisfy two conditions: its methods of proof should be reduced to a minimum number of axioms and rules of inference, and these axioms need to be properly justified. The analysis of the concept of computability made the first condition much more precise. It demands a logic with a set of axioms and rules of inference from which a recursively enumerable set of theorems can be generated (see \cite{godel3}, p. 45). Since the set theory is a demonstrably incomplete theory, the second condition allows that its axioms be given either {\it semantic} or {\it syntactic} justification (see \cite{milos}, pp. 163-168).

Semantic justification is based on ``the insight into the nature of mathematical objects or the understanding of the primitive concepts described by the axioms''. It can be reduced to the verification of set-theoretic axioms in the set universe. Syntactic justification of an axiom is based on its satisfactory consequences. According to G\"{o}del, we are allowed to accept an axiom with a strong syntactic justification, but we should hope to eventually realize that it is semantically justified as well. All $ZF$ axioms have semantic justification, but they do not provide a complete description of the universe of sets. G\"{o}del thought that these axioms are insufficient, not only because they do not allow the proof of $ZF$'s consistency, i.e. the sentence $Con(ZF)$, but also because they cannot decide some mathematically important and natural questions, such as {\it the Continuum Hypothesis} or {\it the Axiom of Choice}. Because of that, they need to be supplemented with some other, properly justified axioms. 

\subsubsection{The Axiom of Constructibility}

G\"{o}del formulates and examines the consequences of one syntactically justified axiom, the so-called {\it Axiom of Constructibility}. In place of the family of sets $(V_{\alpha}: \alpha \in Ord)$, G\"{o}del considers the family of {\it constructible sets} $(L_{\alpha}: \alpha \in Ord)$ whose levels $L_{\alpha}$ are formed either by the power set operation restricted to sets definable by the $ZF$ formulas with parameters and quantifiers pertaining to the lower levels of the hierarchy, or by the union of all $L_{\beta}$, for $\beta < \alpha$ if $\alpha$ is a limit ordinal. It might appear that, by introducing the condition of definability, G\"{o}del gives up on the iterative conception of sets. This is not the case since what he considers is a {\it constructible iterative conception of sets}. If $L$ is the union of the family of constructible sets, the Axiom of Constructibility claims that $V=L$, i.e. that every set is constructible. 

As previously said, $x\in V$ is a $ZF$ formula. If this also holds for $x\in L$, then $V=L$ is a $ZF$ formula as well. If $x\in L$ holds, then there is an ordinal $\alpha$ such that $x\in L_{\alpha}$. The set $x$ consists of the elements $y\in L_{\beta}$ defined by a formula with parameters from $L_{\beta}$, for some $\beta<\alpha$. G\"{o}del has shown that every such $y\in L_{\beta}$ can be obtained by the composition of some of the eleven {\it elementary set operations} applied to the transitive set $L_{\beta}$ and the parameters of defining formula (see \cite{jeh}, pp. 177-178). The universe $L$ is a minimal transitive model of $ZF$ which contains all ordinals and is constructed by an iterative procedure. After the addition of the axiom $V=L$ to $ZF$ theory, its set of axioms remains recursive.

The Axiom of Constructibility has undeniable syntactic justification since its consequences are numerous and include the Generalized Continuum Hypothesis, as well as the Axiom of Choice. After presenting his results at Princeton in 1938, G\"{o}del concludes that it ``seems to give a natural completion of the axioms of set theory, in so far as it determines the vague notion of an arbitrary infinite set in a definite way'' (see \cite{godel2}, p. 27). Besides syntactic justification provided by its important consequences, the axiom $V=L$ has a semantic justification as well, insofar as it is based on the iterative construction of the universe $L$. But, even though he believed in some kind of equilibrium between the semantic and syntactic justification, G\"{o}del still did not find this sufficient for accepting the Axiom of Constructibility. What prevailed seems to be his impression that the constructible universe does not provide a sufficient clarification of the continuum problem, which makes the strongest link between the set theory and the mathematics necessary for understanding the phenomena of the physical world.

The Axiom of Constructibility allows for the existence of inaccessible cardinals, but not of some significantly larger ones. Assuming that there are measurable cardinals, it can be shown that there is a non-constructible set (see \cite{scot}). Cohen's proof of the independence of the Continuum Hypothesis suggested that the semantic justification of the Axiom of Constructibility is insufficient. It seemed that this put an end to the objectivist investigation of the continuum problem, which led the majority of set theorists to formalistic speculations. The interest in a semantic, that is, objective justification of set-theoretic assumptions as understood by G\"{o}del was completely lost. On the other hand, in the letter to Church, G\"{o}del emphasizes his disagreement with the opinion that the proof of its independence from ZFC axioms solves the continuum problem (\cite{godel4}, p. 372). He sticks to his opinion that mathematical concepts and axioms ``describe some well-determined reality'' in which ``Cantor's conjecture must be either true or false, and its undecidability from the axioms as known today can only mean that these axioms do not contain a complete description of this reality'' (see \cite{godel2}, p. 181). 

In the commentary on Cohen's result from 1964, G\"{o}del formulates the program of investigation of large cardinals, which could shed some light on the ``open arithmetic and set-theoretic problems, including the Continuum Hypothesis'' (see \cite{godel3}, pp. 260-261). It is not clear whether, besides the $ZF$ axioms, there are other set-theoretic axioms with semantic justification (see \cite{milos}). But G\"{o}del seems to have thought that such axioms belong to those stating the existence of large cardinals.

G\"{o}del's investigations in set theory seem to be based on the following assumptions: the theory has objective content; its methods are finite and limited by the incomplete formal system; the assumptions about sets need to be undeniably confirmed by our intuition of sets. Simply put, G\"{o}del relies on the assumption of the objectivity of sets and {\it intuition} as a tool that provides insight into their nature within the limits put by reason. Despite the immeasurable significance of his investigations for the development of set theory, G\"{o}del's views on the objectivity of the world of mathematical objects and intuition as a tool that gives insight into its properties were still rejected due to the prevailing materialistic prejudices. But no matter how strong these prejudices are, for mathematicians that accept the existence of the set of natural numbers, which is the weakest form of Platonism, $ZF$ theory has semantically justified axioms, which means that one huge part of mathematics does not have a speculative character at all. 

\section{Rational perception}
On the occasion of his acceptance into the American Philosophical Society in 1961, G\"{o}del was invited to give a talk on a topic of his choice. He prepared a manuscript in which he presented ``in the light of philosophy'' his view on ``the modern development of the foundations of mathematics'' (see \cite{godel3}, pp. 374-387). He begins by noticing that since Renaissance, philosophy was moving away from metaphysics, and closer to positivism, materialism, and skepticism. On the other hand, G\"{o}del does not hide his sympathy towards theology and idealism that ``see sense, purpose and reason in everything'' (see \cite{godel3}, p. 375). He had the opportunity to become very familiar with the principles on which the main trends of positivistic science and philosophy were based since he was a regular participant in the meetings of Vienna Circle, a group of influential mathematicians, logicians, and philosophers whose primary goal was to ``eliminate every form of metaphysics'' from science and philosophy. He never accepted their views and continued building his own inspired by the metaphysical assumptions of mathematics. Contrary to philosophy which was adjusting to the spirit of the time and moving away from the sphere of values, mathematics, as an a priori science was moving in the opposite direction by ``lifting the level of abstractness'' and coming closer to metaphysics. Still, this did not separate it from science, whose connection to mathematics led to a partial change of paradigm, in which ``the possibility of knowledge of the objectivizable states of affairs is denied, and it is asserted that we must be content to predict results of observation'' (see \cite{godel3}, p. 377). 

\subsection{G\"{o}del's rationalism}

At every stage of his professional development and in each of his works, published or unpublished, rationalism appears as G\"{o}del's principal standpoint. It reveals his faith in human reason and belief that it can answer in a systematic way all meaningful questions it poses. Although strict, his rationalism is much more encompassing than usually understood. It pertains not only to science and philosophy but to the entirety of our knowledge. In a conversation with Carnap from 1940, G\"{o}del says that ``One could establish an exact system of postulates employing concepts that are usually considered metaphysical: `God', `soul', `idea'. If this is done accurately, there would be no objection'' (see \cite{gierer} for Carnap's notes of the conversation). Responding to Carnap's remark that this is the case if the system is understood as ``just a calculus'', G\"{o}del says explicitly that he is ``thinking not of a calculus, but of a theory... [which] has observable consequences, but the observed consequences do not exhaust the theory''. He remarks in the same conversation that ``decisive progress in science - including physics - often depends on a change of direction''. In analogy to Leibniz's {\it Calculus Ratiocinator}, G\"{o}del believes that it is possible to build a deductive system encompassing the totality of our knowledge, not only scientific but also that belonging to the area of values and metaphysical ideas. According to him, positivistic science and metaphysics are not mutually exclusive and there is no need for their sharp separation. On the contrary, their unified development would be profitable for both. 

G\"{o}del's rationalism naturally applies to mathematics as well. In that form, it resembles Hilbert's conviction that ``every definite mathematical problem must necessarily be susceptible of an exact settlement, either in the form of an actual answer to the question asked, or by the proof of the impossibility of its solution'' (\cite{hilbert1}, p. 444). This conviction led Hilbert to formulate his program of founding the whole of mathematics on a finite set of arithmetic axioms whose consistency was supposed to be provable from them alone. Although incompleteness theorems show that this cannot be done, in G\"{o}del's opinion, they do not show that Hilbert's rationalism is unfounded. ``If my result is taken together with the rationalistic attitude which Hilbert had and which was not refuted by my results, then [we can infer] the sharp result that mind is not mechanical'' (\cite{wang1}, 6.1.8). Instead of refuting rationalism in mathematics, incompleteness theorems led G\"{o}del to a particular philosophical view on the way mathematical knowledge is attained and extended. 

\subsection{Conceptual realism and its epistemology}

While G\"{o}del's rationalistic attitude is clear and indubitable, his Platonism makes for a much more complex position. He develops this standpoint very carefully, relying on the truth of mathematical propositions, syntactic or semantic, which would not be possible if the terms appearing in them and determining their truth value lacked reference. G\"{o}del's philosophy is best characterized by his {\it conceptual realism}. As we said, G\"{o}del was always a Platonist, but his Platonism at first had the form of mathematical realism which, as the final result of his philosophical contemplation, developed into conceptual realism. Together with it, G\"{o}del formed an opinion on the nature of mathematical knowledge as ultimately based on understanding objective mathematical concepts. His understanding of the nature of abstract ideas is not confined to the mathematical realm but applies also to the mathematical descriptions of physical reality: ``Mathematical propositions, it is true, do not express physical properties of the structures concerned, but rather properties of the {\it concepts} in which we describe those structures. But this only shows that the properties of those concepts are something quite as objective and independent of our choice as physical properties of matter'' (\cite{godel3}, p. 360). With the hope of making his position stronger, G\"{o}del was trying for quite a long time to establish an adequate epistemology of conceptual realism. Since the formal tools do not provide a complete insight into the mathematical world, G\"{o}del assumed that they cannot be the only way to gain knowledge of it. In addition, we need to rely upon some kind of direct insight into the content of the corresponding concepts. We need some kind of perception analogous to the direct sensory experience of the material world. G\"{o}del believed that the mind is capable of such a perception of the conceptual world. His confidence in the existence of this kind of mental ability is a consequence of his belief that our knowledge is not limited by the mechanical nature of the formal systems but can grow ``in every direction''. 

G\"{o}del's most comprehensive work on the consequences of incompleteness theorems is the manuscript of his lecture from 1951 (see \cite{godel3}, pp. 304-323). One part of this lecture is devoted to the philosophy of mind in which G\"{o}del presents the consequences of incompleteness theorems in the form of dichotomy: {\it either the human mind infinitely surpasses the powers of any finite machine, or else there exist absolutely unsolvable diophantine problems}. Absolutely unsolvable problems are unsolvable not just in a particular formal system, but with the help of ``any imaginable mathematical decision proof''. It is worth noting that the dichotomy is not exclusive. G\"{o}del characterizes this statement as a ``mathematically established fact'' and he concludes that both alternatives have consequences that are ``definitely opposed to materialistic philosophy''. The first alternative implies that thinking cannot be reduced to a formal system, while the second one refutes the view that mathematics is a creation of some ideal mathematical subject. ``Creator necessarily knows all properties of his creatures ... so this alternative seems to imply that mathematical objects and facts (or at least {\it something} in them) exist objectively and independently of our mental acts and decisions'' (\cite{godel3}, p. 311). Even though it speaks in favor of the objectivist standpoint, G\"{o}del rejects the alternative that implies the existence of absolutely unsolvable mathematical problems. He had very different reasons for believing in mathematical Platonism and was inclined to accept the first alternative. Still, he was careful since he knew that both mechanicism and the existence of absolutely unsolvable problems are consistent with his incompleteness theorems. His main reasons for rejecting the second alternative are mostly philosophical. He accepts Kant's optimistic standpoint, which implies that our mind would be ``fatally irrational'' if it would be capable of posing the questions to which it cannot answer. By accepting the first alternative, G\"{o}del justifies his belief in the existence of mathematical {\it intuition} that surpasses the cognitive capacities of formal systems. 

We should bear in mind that G\"{o}del never presented his philosophical views systematically. Instead, they are scattered throughout his works and often formulated as remarks in different contexts, which makes the credibility of their interpretation questionable. Some of them are personal notes never planned to be published. This applies primarily to those pertaining to G\"{o}del's attempts at founding the epistemology of the objectivist understanding of the world of concepts. Although, as we said, rationalism and Platonism are the constants in G\"{o}del's philosophy, this cannot be said for his epistemological views. His understanding of intuition was changing according to developments in logic and set theory, as well as the teachings of other philosophers, and it never assumed a definite shape. It is only indicated in a number of informal parallels between intuition and sense perception, and intuition and reason, where intuition is often taken to have the role in the discovery of set-theoretic axioms that should decide the open problems of this theory. Even though G\"{o}del was using these concise analogies only to problematize the idea of intuition, in the philosophical literature, they received a detailed analysis and refutation as if they expressed G\"{o}del's formed position (see \cite{parsons}). This applies, for example, to G\"{o}del's following remark: ``we use senses to perceive particular objects, (while) with the mind we perceive concepts (in particular, primitive ones) and their relations''. Since the sets are the way of ``understanding concepts'', when he discusses intuition, G\"{o}del primarily has in mind the intuition of sets. ``Despite their remoteness from sense experience, we do have something like a perception also of the objects of set theory, as is seen from the fact that the axioms force themselves upon us as being true. I don't see any reason why we should have less confidence in this kind of perception, i.e., in mathematical intuition, than in sense perception...'' (see \cite{godel2}, p. 268). These remarks about the {\it possibility} of some kind of mental ability that allows us to perceive mathematical objects are nothing but G\"{o}del's attempts at supplementing his conceptual realism with the appropriate epistemology. He tried to solve this problem by relying on Kant's views on intuition, and later, on Husserl's phenomenology, in which he put a lot of hope, but which turned out not to be useful in forming his own views.

But initially, G\"{o}del seems to have been relying on the pre-Kantian understanding of intuition, which does not clearly separate the philosophy of mind and epistemology, since he sees intuition as a cognitive tool the mind uses to surpass the power of a formal (mechanical) system. At this period, his understanding of intuition was most closely related to Leibniz's view according to which knowledge is intuitive if it is:
\begin{itemize}
    \item {\it clear}, i.e. it gives the means for recognizing the object it concerns;
    \item {\it distinct}, i.e. one is in a position to enumerate the marks or features that distinguish an instance of one’s concept;
    \item {\it adequate}, i.e. one’s concept is completely analyzed down to primitives;
\end{itemize}
and finally, if all these elements are understood immediately and momentarily (see \cite{lajbnic}, p. 23). Lebniz's definition tells us that a complete, rational analysis of an object is supposed to provide us with a kind of apprehension of its characteristics. Having in mind this kind of rational analysis, G\"{o}del at first relates intuition to a kind of apprehension and bases his analogy with sense perception on it. However, later he insists that ``this additional sense (i.e., reason) is not counted as a sense, because its objects are quite different from those of all other senses. For, while through sense perception we know particular objects and their properties and relations, with mathematical reason we perceive the most general (namely the ``formal'') concepts and their relations, which are separated from space-time reality'' (see \cite{godel3}, p. 354). G\"{o}del replaces the idea of intuition as some kind of apprehension that provides us with new knowledge with {\it rational perception}. Just as we gain knowledge in set theory by first developing it, and then adding some new axioms for which we seek semantic and syntactic justification, so too is the knowledge in other mathematical fields gained by working out the properties of primitive concepts and thus explicating their content, and by formulating new assumptions that should decide the problems that remained unsolved. Because of the inexhaustibility of mathematical knowledge, the number of these rational perceptions is unlimited (see \cite{godel3}, p. 353). G\"{o}del built his view of intuition as rational perception taking as a paradigm the investigations in set theory that consist of the introduction of new axioms, their semantic justification in the iterative conception of set, and the verification provided by the consequences they have in set theory. 

\section{Alternative philosophies of mathematics}

Directly opposed to G\"{o}del's understanding of the objectivity of mathematics is {\it nominalistic standpoint}, which denies mathematical propositions and knowledge any objective content. The main nominalistic argument against Platonism is the {\it assumption of causality} in the definition of knowledge, which is absent from the Plato's original definition according to which knowledge is {\it justified true belief}. Since mathematical objects are acausal, i.e. not belonging to space or time, nominalists deny that mathematical knowledge is knowledge at all. On the other hand, given the status of mathematical knowledge as indisputable and the foundation of every serious study of the phenomena in the world, it would be rather natural that the assumption of causality is refuted and some other, more appropriate reformulation of Plato's definition is asked for. The way of justifying a true belief differs across the areas of knowledge. Thus, in mathematics, a true belief is justified by its {\it proof}, in science by (experimental) {\it verification}, and in philosophy by {\it argumentation}. It is clear that proof provides a more reliable confirmation of knowledge than scientific verification or philosophic argumentation. Ever since antiquity, all other areas aspired to reach the exactness that characterizes mathematical knowledge, i.e. to verify their knowledge in the way mathematics does. It is difficult to imagine what sciences would look like if they hadn't been closely related to mathematics. They are so dependent on mathematics that it is often difficult to discern the line that separates them from it. It is exactly the acausal, or, as G\"{o}del would say, ``analytic nature'' that grants precision to the mathematical models of spatially and temporally determined phenomena. 

Different philosophies of mathematics put restrictions on the content of mathematical propositions, motivated by particular epistemological scruples. The real content of mathematical propositions is either denied (by nominalism, formalism, or conventionalism), or is only partially admitted, while any use of nonconstructive arguments, such as the law of excluded middle, the Axiom of Choice, or the existence of some infinite sets, is deemed illegitimate (by intuitionism, constructivism, or restricted realism of Quinean type). G\"{o}del had the idea to provide a strong defense and justification for his objectivist conception of mathematics, i.e. for conceptual realism, by: (1) making a list of all the alternative philosophies of mathematics, (2) proving that the list is exhaustive, and (3) showing that his incompleteness theorems rule out every option from it. Namely, each of these conceptions has the ambition to ground the whole of mathematics. Mathematics would thus have to be expressed in a formalism that contains arithmetic, which would make it incomplete. This argument helped G\"{o}del to show that intuitionism, conventionalism, and other views are untenable. However, he did not manage to show the completeness of the list of alternative philosophies. Being aware of this, G\"{o}del still argued that this argument can be conducted ``with mathematical rigor'' (see \cite{godel3}, p. 322). The weakness of G\"{o}del's argument could be overcome by dividing all options from the list into two classes, e.g. into the nominalistic and realistic philosophies. It could then be shown that incompleteness theorems threaten all the options from these two classes as long as they do not allow for any other way of arriving at mathematical knowledge besides the one contained in its formalization. Mathematical realism would not be exempt from this criticism. But G\"{o}del's conceptual realism is supposed to overcome this. The implications of this particular philosophical standpoint are supposed to show in a future theory that G\"{o}del called the logic of concepts. If we managed to establish this theory successfully and if it turned out to be revealing about mathematical subjects as well, this would make for a strong argument in favor of G\"{o}del's view that mathematical knowledge is conceptual in nature. Before we say something about how this theory is supposed to look like, we will review what in our opinion are G\"{o}del's primary reasons for insisting on its development.

\section{Incompleteness and extensionality}

From what is said in the previous chapters it follows that there are at least two reasons for incompleteness of set theory. The first one lies in its subject. It deals with extensional objects, that is, sets taken as simple collections with no structure or a rule according to which they are built. It is difficult to obtain a complete description of the totality of such objects. The second reason lies in the way set theory deals with this subject, and it affects all other theories of the same kind. Namely, it is a first-order theory with a recursive set of axioms containing arithmetic. As G\"{o}del has shown, every such theory, if consistent, will be incomplete and its consistency will not be provable in it. In what follows, we will analyze both reasons more thoroughly to determine if and how they could be overcome. As we will try to show, it is the focus on objects and extensions (instead of concepts) that makes completeness unattainable. 

Let us first consider the subject of set theory. The interpretation of $ZF$ theory provides a specific view of sets and their formation. According to this interpretation, sets constitute a cumulative hierarchy that respects the order in which they are formed. Since a set can be formed only after all of its elements, the first level of the hierarchy can only be occupied by the empty set (or some objects that are not sets). For the same reason, sets from the higher levels can contain only sets from the lower levels as their elements. Also, every newly formed set gives rise to some other sets containing it. This process of set formation can be continued indefinitely. Consequently, the hierarchy of sets can indefinitely grow.

$ZFC$ axioms describe only a part of this hierarchy. This is made obvious by the existence of propositions stating some properties of sets belonging to this hierarchy (such as {\it CH}), which are undecidable from the axioms. The description provided by set theory can be amended with some new axioms. As a rule, they state the existence of different large cardinals. They thus describe the properties of some very high levels of the hierarchy that cannot be inferred from what is known about the sets in the lower levels. On the other hand, these axioms can be quite revealing about the properties of the smaller sets.

However, we cannot hope that the additions of such axioms will eventually lead to a complete theory. The reason is that set theory describes the building of the hierarchy of sets step-by-step. This follows from the bottom-up, inductive building of sets described by the axioms based on the iterative conception (cf. \cite{godel3}, pp. 306-307). Such description, motivated by the thoroughly extensional understanding of sets, does not provide us with a way of inferring the properties of all sets, and it does not allow us to understand the set hierarchy in its entirety. So it cannot be hoped that it will eventually yield a theory that cannot be extended by axioms that decide some new propositions about sets (cf. \cite{godel3}, pp. 45-48). 

In his Princeton lecture, G\"{o}del remarks that one way of reaching completeness of this kind would be finding a nonconstructive characterization of the large cardinal axioms whose addition to set theory has significant consequences. In this way, all the steps of building the hierarchy are supposed to be described in a nonconstructive way. The notion of proof from the existing and any of these new axioms would then make for an approximation of the notion of proof in a complete set theory (\cite{godel2}, p. 151). But how to arrive at this characterization? According to G\"{o}del, the axioms of set theory are supposed to describe the properties of {\it the concept of set}. Every new axiom can thus be seen as a step in the development of our understanding of this concept. The description of all the axioms by which set theory can be extended would then be based on our complete understanding of the concept of set. So the decisive step from incomplete to complete theories would seem to essentially involve the transition from extensional to intensional considerations. As G\"{o}del remarks: ``A complete foundation of set theory calls for a study of properties and concepts'' (\cite{wang1}, 8.6.22). This study should lead to a better understanding of the concept of set, which should be irreducible to its extensional description. Namely, what the incompleteness of mathematical formal systems might be taken to show is that if we try to translate all the knowledge we gain by understanding a particular concept into the knowledge of its extension, we are bound to incompleteness. 

Taking into account some nonconsecutive characterization of axioms also means changing the methods of set theory. This could lead to the abolishment of the second reason for its incompleteness, according to which set theory is incomplete because it is a sufficiently strong first-order theory with a recursive set of axioms. The explanation for the fact that all such theories are incomplete, if they are consistent, might again be that they are dealing with particular mathematical concepts by focusing only on their extensions. First-order theories give us methods for dealing with objects. We can deal with properties and concepts only insofar as we can reduce them to objects or their collections. They allow us to systematically derive all the consequences of our basic insights formulated in axioms, which, although reached by understanding related concepts, concern the objects to which these concepts apply. But, they do not allow us to go beyond them. On the other hand, our knowledge is capable of developing since ``we understand abstract terms more and more precisely as we go on using them, and ... more and more abstract terms enter the sphere of our understanding'' (\cite{godel2}, p. 306). There is, thus, another way of reaching knowledge, which is different from any mechanical procedure. It involves the cultivation of our understanding of concepts. 

The transition to these intensional considerations would thus involve methods mathematicians use to develop their understanding of the subject matter, which cannot be imitated by any machine, but might include what G\"{o}del used to call ``mathematical intuition'' or ``rational perception''. According to G\"{o}del, this process ``today is far from being sufficiently understood to form a well-defined procedure. It must be admitted that the construction of a well-defined procedure which could actually be carried out (and would yield a non-recursive number-theoretic function) would require a substantial advance in our understanding of the basic concepts of mathematics'' (\cite{godel2}, p. 306). G\"{o}del might have hoped for this to be achievable in a theory that describes the formal properties of concepts and their mutual relations. It could show which aspects of concepts as well as the relations between them are crucial for their adequate and complete understanding. Consequently, it might help us determine the limits of the formal treatment of mathematical concepts. This might be part of the reason why G\"{o}del insisted on the importance of founding such a theory.

In the next section, we review G\"{o}del's suggestions about how the logic of concepts should look like. We try to determine if some kind of completeness could be achievable in it, and if that would have any implications for the problem of completeness of set theory.

\section{The logic of concepts}

G\"{o}del envisioned the logic of concepts as a theory that deals with the formal properties of concepts understood as intensional meanings of predicates. They should thus contain all and only the properties an object needs to have for this predicate to be truthfully attributed to it (i.e. to its name). The basic relation this theory is supposed to characterize is the relation of {\it application}. This is the relation in which a concept stands to objects that satisfy properties contained in it. It should be essentially different from the membership relation in which the elements stand to a set. According to G\"{o}del, their difference is best revealed by the fact that the relation of concept application is not irreflexive, which means that a concept can apply to itself. Examples of such concepts are {\it the concept of concept}, {\it the concept of concept with infinite range}, {\it the concept of abstract thing}, etc (cf. \cite{wang1}, 8.6.3). This self-applicability of concepts is possible because of their intensional nature, that is, the fact that they are not built out of the objects to which they apply (in contrast to sets that are built out of their elements). The objects to which a concept applies do not uniquely distinguish this concept, since two different concepts may apply to the same things. What exactly distinguishes them, and what are the criteria of their identity, is a question to be assessed in the logic of concepts. 

The self-applicability of concepts can lead to paradoxes similar to those that appeared in the so-called na\"{i}ve set theory. For example, the intensional version of Russell's paradox concerns the self-applicability of {\it the concept of not applying to itself}. According to G\"{o}del, these paradoxes show that concepts are objective since we are not free to construct and apply them as we wish (see \cite{wang1}, 8.5.20). If we hope to establish a theory dealing with them, we need first to find a solution to the paradoxes that is based on the nature of concepts. It is thus supposed to reveal some important properties of concepts and their application that justify particular restrictions on their formation. 

One way of dealing with this problem would be to treat the relation of concept application as a {\it partial} relation. This implies allowing for the possibility that neither a concept nor its complement applies to the given object. In other words, the question of whether a concept applies to an object or not, would not always be taken as meaningful. Contrary to that, the question of whether a set contains an object or not is always meaningful. 

This way of dealing with intensional paradoxes is pointed out by some of G\"{o}del's remarks (see \cite{godel2}, 137-138). They suggest that we should try to solve intensional paradoxes by scrutinizing the notion of meaningfulness, trying to find out what makes an application of a concept meaningful or not, and if there is a criterion of meaningfulness that rules out the applications of concepts that lead to paradoxes, but is not so restrictive to lead to the theory of types.

The partiality of the relation of concept application makes concepts similar to computable functions. It is well-known that some computable functions do not yield any result when applied to particular arguments. In other words, the process of computation does not always terminate. Also, it cannot always be determined whether the computation terminates or not. Similarly, it might not always be possible to say if a concept applies meaningfully to an object. This accords well with the intensionality of both concepts and computable functions. It shows that they are not determined by the objects to which they apply (and the results this yields), but rather by the rules of this application. 

If a concept does not always meaningfully apply to an object, then its {\it range}, formed by objects to which it applies, and {\it anti-range}, formed by objects to which its complement applies, need to be treated separately. The axioms of the logic of concepts would then be supposed to describe the relation in which a concept stands to the objects from its range, as well as to those from its anti-range (see \cite{kostic} for some suggestions on how this can be done). The theory is also supposed to describe how complex concepts are built out of primitive ones, using logical connectives. This description might lead to establishing a criterion of identity of complex concepts, which might be reducible to an equivalence relation between their structures determined by this formation, or between the primitive concepts from which they are built.

Given the subject of the envisioned theory and the questions it is supposed to answer, we might wonder if it could be built as a complete theory, and what would this completeness stand for.

\subsection{Intensionality and completeness}

If the incompleteness of set theory has to do with the extensionality of its subject and method, as suggested by the above-given explanation, then it might be avoidable in a theory that deals with an intensional subject. The logic of concepts would presumably strive for completeness (see \cite{kostamilos}, p. 43), whose achievement might have some significant consequences for set theory as well, at least under the assumption of G\"{o}del's mathematical epistemology. 

To encourage the quest for the absolute definition of provability, in which the logic of concepts might have an important role, G\"{o}del gives the example of computability theory that manages to provide an absolute definition of the related concept of computability. In analogy to computability theory, a complete theory of concepts is supposed to contain all the available principles for deciding if a concept built in a particular way is meaningfully applicable to some object or not. In other words, it would not be possible to extend it by some additional principles that would decide more questions about the meaningful applicability of concepts. The reason why something similar is not possible in set theory here does not hold, because the logic of concepts is not supposed to deal with a process that cannot be finished, nor described in a way that allows us to infer all the properties of the objects resulting from it. Since it is supposed to deal with intensional objects, we may expect it to establish some principles that would allow us to understand the properties of these objects, which do not change with the growth of their range or the complexity of their construction. 

Given G\"{o}del's epistemic optimism, he might have hoped for still stronger completeness of the logic of concepts. That is, he might have expected that it would be able to prove every true proposition about its subject. For this to be possible, the methods of the logic of concepts would have to differ from those of set theory. G\"{o}del seems to have imagined the logic of concepts as a theory formulated in the predicate calculus by adding the relation of application (\cite{wang1}, 8.6.18). If this theory is supposed to contain arithmetic, then the only way for G\"{o}del's incompleteness theorems not to apply to it would be that its set of axioms is non-recursive (this accords well with G\"{o}del's remark that a well-defined procedure for understanding concepts would yield a non-recursive number-theoretic function). That something is an axiom of this theory would not then always be decidable by examining its formal structure. The meaning of the terms appearing in them and our understanding of that meaning, which cannot be imitated by a machine, would also need to be called for.

\subsection{The logic of concepts and set theory} 

G\"{o}del remarks that the future logic of concepts might contain set theory. This is explained by his conjecture that for every set, there is a (defining) concept (\cite{wang1}, 8.6.4). In that case, we would be able to deal with sets as extensions of the corresponding concepts. According to G\"{o}del, we might have ``something like a hierarchy of concepts (or also of classes) which resembles the hierarchy of sets and contains it as a segment'' (\cite{wang1}, 8.6.20). But he stresses that this hierarchy would be peripheral to the theory of concepts and could not be taken as the appropriate interpretation of this theory. 

There is another, indirect way in which concept theory might improve our knowledge of sets and other mathematical objects. Namely, by developing our understanding of the universal properties of concepts and their relations, the logic of concepts could contribute to the adequate understanding of mathematical concepts (including {\it the concept of set}), and to recognizing all the properties of mathematical objects that follow from it. G\"{o}del remarks that {\it the concept of concept} and {\it the concept of absolute proof} are closely related (cf. \cite{wang1}, 6.1.13). Taking into account his mathematical epistemology, this is not so difficult to understand. It seems to imply that every notion of proof in mathematics is based on some understanding of related concepts. A complete understanding of those concepts should enable us to establish a theory that provides an absolute notion of provability in the given area. This suggests that the logic of concepts is supposed to reveal how we base proofs of mathematical propositions on the content of corresponding concepts. It would thus be concerned, among other things, with deductions tied to concepts (see \cite{kostamilos}, pp. 43-47). In that way, it could improve the methods by which mathematicians reach new knowledge of a mathematical concept and derive its consequences for the objects falling under it. 

\section{Conclusion}

Thanks to G\"{o}del's results, the content and the method of modern logic drastically differ from those that characterized it for millennia. He defined the concept of formal logical system; he precisely distinguished the main problems of such systems - their consistency, completeness, decidability, and categoricity; and he investigated or solved those problems for the formal systems important for the foundations of mathematics. He arrived at these results relying on philosophical assumptions that do not accord with the positivistic and naturalistic spirit of our time, which is why philosophers denied him their support in developing them. Encouraged by the fruitfulness of his assumptions and fully aware of the destructive consequences of positivistic and naturalistic assumptions for mathematics, G\"{o}del persisted in his belief in the objective existence of the mathematical and conceptual world. We believe that philosophy which is not limited by the assumption of the exclusive existence of the perceptible world and linguistic conventions used in its description could contribute much more to understanding G\"{o}del's ideas contained in his unpublished works, such as his idea of the logic of concepts we dealt with in this article.

\bibliographystyle{amsplain}

\end{document}